\documentclass[reqno]{amsart}

\usepackage{amssymb}
\usepackage{graphicx}
\usepackage{amscd}
\usepackage[hidelinks]{hyperref}
\usepackage{color}
\usepackage{float}
\usepackage{graphics,amsmath,amssymb}
\usepackage{amsthm}
\usepackage{amsfonts}
\usepackage{latexsym}
\usepackage{epsf}
\usepackage{xifthen}
\usepackage{mathrsfs}
\usepackage{dsfont}
\usepackage{makecell}
\usepackage[FIGTOPCAP]{subfigure}
\usepackage{amsmath}
\allowdisplaybreaks[4]
\usepackage{listings}
\usepackage{etoolbox}
\usepackage{fancyhdr}
\usepackage{pdflscape}
\usepackage[title,toc,titletoc]{appendix}
\usepackage{enumitem}
\usepackage[noadjust]{cite}

 \newtheoremstyle{mytheorem}
 {3pt}
 {3pt}
 {\slshape}
 {}
 {\bfseries}
 {.}
 { }
 {}

\numberwithin{equation}{section}

\theoremstyle{theorem}
\newtheorem{theorem}{Theorem}[section]
\newtheorem*{theorem*}{Theorem}
\newtheorem{corollary}[theorem]{Corollary}
\newtheorem{lemma}[theorem]{Lemma}

\theoremstyle{definition}

\newtheorem*{example*}{Example}
\newtheorem{conjecture}{Conjecture}[section]

\theoremstyle{remark}

\newtheorem*{remark*}{Remark}
\newtheorem*{remarks*}{Remarks}

\newcommand{\doi}[1]{\href{http://dx.doi.org/#1}{doi: #1}}

\newcommand{\Keywords}[1]{\ifthenelse{\isempty{#1}}{}{\smallskip \smallskip \noindent \textbf{Keywords}. #1}}
\newcommand{\MSC}[2][2010]{\ifthenelse{\isempty{#2}}{}{\smallskip \smallskip \noindent \textbf{#1MSC}. #2}}
\newcommand{\abstractnote}[1]{\ifthenelse{\isempty{#1}}{}{\smallskip \smallskip \noindent \textsuperscript{\dag}#1}}

\makeatletter
\def\specialsection{\@startsection{section}{1}%
  \z@{\linespacing\@plus\linespacing}{.5\linespacing}%
  {\normalfont}}
\def\section{\@startsection{section}{1}%
  \z@{.7\linespacing\@plus\linespacing}{.5\linespacing}%
  {\normalfont\scshape}}
\patchcmd{\@settitle}{\uppercasenonmath\@title}{\Large\boldmath}{}{}
\patchcmd{\@settitle}{\begin{center}}{\begin{flushleft}}{}{}
\patchcmd{\@settitle}{\end{center}}{\end{flushleft}}{}{}
\patchcmd{\@setauthors}{\MakeUppercase}{\normalsize}{}{}
\patchcmd{\@setauthors}{\centering}{\raggedright}{}{}
\patchcmd{\section}{\scshape}{\large\bfseries\boldmath}{}{}
\patchcmd{\subsection}{\bfseries}{\bfseries\boldmath}{}{}
\renewcommand{\@secnumfont}{\bfseries}
\patchcmd{\@startsection}{\@afterindenttrue}{\@afterindentfalse}{}{}
\patchcmd{\abstract}{\leftmargin3pc}{\leftmargin1pc}{}{}

\def\maketitle{\par
  \@topnum\z@ 
  \@setcopyright
  \thispagestyle{empty}
  \ifx\@empty\shortauthors \let\shortauthors\shorttitle
  \else \andify\shortauthors
  \fi
  \@maketitle@hook
  \begingroup
  \@maketitle
  \toks@\@xp{\shortauthors}\@temptokena\@xp{\shorttitle}%
  \toks4{\def\\{ \ignorespaces}}
  \edef\@tempa{%
    \@nx\markboth{\the\toks4
      \@nx\MakeUppercase{\the\toks@}}{\the\@temptokena}}%
  \@tempa
  \endgroup
  \c@footnote\z@
  \@cleartopmattertags
}
\makeatother

\newcommand{\bI}{\mathbf{I}}
\newcommand{\cD}{\mathcal{D}}
\newcommand{\cF}{\mathcal{F}}
\newcommand{\cL}{\mathcal{L}}
\newcommand{\cR}{\mathcal{R}}
\newcommand{\last}{\textsc{last}}
\newcommand{\srpt}{\textsc{srpt}}

\title[On $0012$-avoiding inversion sequences]{On $0012$-avoiding inversion sequences and a Conjecture of Lin and Ma}

\author[S. Chern]{Shane Chern}
\address{Department of Mathematics, Penn State University, University Park, PA 16802, USA}
\email{shanechern@psu.edu; chenxiaohang92@gmail.com}

\date{}

\begin{document}

%

\maketitle

\begin{abstract}

The study of pattern avoidance in inversion sequences recently attracts extensive research interests. In particular, Zhicong Lin and Jun Ma conjectured a formula that counts the number of inversion sequences avoiding the pattern $0012$. We will not only confirm this conjecture but also give a formula that enumerates the number of $0012$-avoiding inversion sequences in which the last entry equals $n-1$.


\MSC{05A05, 05A15.}
\end{abstract}

\section{Introduction}

An \textit{inversion sequence} of length $n$ is a sequence $e=e_1 e_2 \cdots e_n$ such that $0\le e_i \le i-1$ for each $1\le i\le n$. We denote by $\bI_n$ the set of inversion sequences of length $n$. Given any word $w \in \{0,1,\ldots,n-1\}^n$ of length $n$, we define its \textit{reduction} by the word obtained via replacing the $k$-th smallest entries of $e$ with $k-1$. For instance, the reduction of $0023252$ is $0012131$. We say that an inversion sequence $e$ \textit{contains} a given pattern $p$ if there exists a subsequence of $e$ such that its reduction is the same as $p$; otherwise, we say that $e$ \textit{avoids} the pattern $p$. For instance, $0023252$ has a subsequence $022$ whose reduction is $011$ --- hence, $0023252$ contains the pattern $011$. On the other hand, none of the length $3$ subsequences of $0023252$ have the reduction $110$ --- hence, $0023252$ avoids the pattern $110$.

Let $p_1$, $p_2$, \ldots, $p_m$ be given patterns. We denote by $\bI_n(p_1,p_2,\ldots,p_m)$ the set of inversion sequences of length $n$ that avoid all of the patterns $p_1$, $p_2$, \ldots, $p_m$. Recently, the study of pattern avoidance in inversion sequences attracts extensive research interests. See \cite{Aul2020,AE2019a,AE2019b,AE2020,BBGR2019,CJL2019,CMSW2016,KL2017,Lin2018,Lin2020,LF2020,LY2020,MS2015,MS2018,YL2019,Yan2018} for several instances of work on this topic. Among these work, one particular interesting problem is about the enumeration of inversion sequences that avoids fixed patterns. For example, in a pioneering work of Corteel, Martinez, Savage and Weselcouch \cite{CMSW2016}, it was shown that
$$|\bI_n(011)|=B_n \quad\text{and}\quad |\bI_n(021)|=S_n$$
where $B_n$ is the $n$-th Bell number (OEIS, \cite[A000110]{OEIS}) and $S_n$ is the $n$-th large Schr\"oder number (OEIS, \cite[A006318]{OEIS}).

In a recent paper \cite{YL2019}, Yan and Lin proved a conjecture due to Martinez and Savage \cite{MS2018} that claims
\begin{equation}
|\bI_n(021,120)|=1+\sum_{i=1}^{n-1}\binom{2i}{i-1}.
\end{equation}
This sequence is registered as OEIS, \cite[A279561]{OEIS}. Lin and Yan also showed that this sequence as well enumerates $|\bI_n(102,110)|$ and $|\bI_n(102,120)|$. This therefore establishes the Wilf-equivalence
\begin{equation}\label{eq:Wilf-equiv}
\bI_n(021,120)\sim \bI_n(102,110) \sim \bI_n(102,120).
\end{equation}
At the end of \cite{YL2019}, a conjecture of Zhicong Lin and Jun Ma discovered in 2019 is recorded.

\begin{conjecture}[Lin and Ma]\label{conj:Lin-Ma}
	For $n\ge 1$,
	\begin{equation}
	|\bI_n(0012)| = 1+\sum_{i=1}^{n-1}\binom{2i}{i-1}.
	\end{equation}
\end{conjecture}

In other words, it is possible to extend the Wilf-equivalence \eqref{eq:Wilf-equiv} as
$$\bI_n(0012)\sim \bI_n(021,120)\sim \bI_n(102,110) \sim \bI_n(102,120).$$

In this paper, we will prove the above conjecture of Lin and Ma.

\begin{theorem}\label{th:0012}
	Conjecture \ref{conj:Lin-Ma} is true.
\end{theorem}

Let us fix some notation. Given $e=e_1e_2\cdots e_n\in \bI_n(0012)$, we define
$$\cR(e):=\{m: \exists\, i\ne j \text{ such that }e_i=e_j=m\}.$$
In other words, $\cR(e)$ is the set of letters that appear more than once in $e$. We further define
$$\srpt(e):=\min \cR(e),$$
that is, the smallest number in $\cR(e)$. Notice that there is only one sequence $01\cdots (n-1)$ in which none of the letters repeat. For this sequence, we assign that
$$\srpt(01\cdots (n-1)):=n-1.$$
Finally, we define
$$\last(e):=e_n,$$
the last entry of $e$.

Apart from counting the number of inversion sequences that avoid the pattern $0012$, we will also enumerate the number of sequences in $\bI_n(0012)$ in which the last entry equals $n-1$.

\begin{theorem}\label{th:0012-last}
	For $n\ge 1$,
	\begin{equation}
	|\{e\in \bI_n(0012): \last(e)=n-1\}|=\begin{cases}
	1 & \text{if $n=1$},\\
	2^{n-2} & \text{if $n\ge 2$}.
	\end{cases}
	\end{equation}
\end{theorem}

\section{Combinatorial observations}

We collect some combinatorial observations about inversion sequences in $\bI_n(0012)$.

\begin{lemma}\label{le:combin-0}
	For $n\ge 1$ and $e\in \bI_n(0012)$, if $\srpt(e)=k$, then for $1\le i\le k+1$,
	$$e_i=i-1.$$
\end{lemma}

\begin{proof}
	If $\srpt(e)=n-1$, then $e=01\cdots (n-1)$ and hence the lemma is true. Let $\srpt(e)\ne n-1$. If in this case the lemma is not true, then since $0\le e_i\le i-1$ for each $i$, there must exist some $k_1<k=\srpt(e)$ that appears more than once among $e_1$, $e_2$, \ldots, $e_{k+1}$. This violates the assumption that $\srpt(e)=k$.
\end{proof}

\begin{lemma}\label{le:combin-1}
	For $n\ge 2$ and $e=e_1e_2\cdots e_n\in \bI_n(0012)$, let $\gamma(e)=e_1e_2\cdots e_{n-1}$. We further assume that $e\ne 01\cdots (n-1)$. Then
	
	\begin{enumerate}[leftmargin=*,align=left,labelsep=0cm,label=\textup{(\alph*)}.]
		\item if $\last(e)>\srpt(\gamma(e))$, then
		$$\srpt(e)=\srpt(\gamma(e));$$
		
		\item if $\last(e)\le \srpt(\gamma(e))$, then
		$$\srpt(e)=\last(e).$$
	\end{enumerate}
\end{lemma}

\begin{proof}
	A simple observation is that $\gamma(e)\in \bI_{n-1}(0012)$. Below let us assume that $\last(e)=\ell$, $\srpt(e)=k$ and $\srpt(\gamma(e))=k'$.
	
	First, if $\cR(\gamma(e))=\emptyset$, then for each $0\le i\le n-1$, $e_i=i-1$. Since $e\ne 01\cdots (n-1)$, we have $\last(e)=\ell\le n-2=\srpt(\gamma(e))$. This fits into Case (b). Further, we find that $\cR(e)=\{\ell\}$ and hence $\srpt(e)=\ell$. This implies that $\srpt(e)=\last(e)$.
	
	Now we assume that $\cR(\gamma(e))\ne \emptyset$. Notice that Case (a) is trivial. For Case (b), we first deduce from $\cR(\gamma(e))\ne \emptyset$ that $k' \le n-3$. By Lemma \ref{le:combin-0}, we find that for $1\le i\le k'+1$, $e_i=i-1$. If $\last(e)=\ell \le k'$, then we know that $e_{\ell+1}=\ell=e_n$. Also, we notice that the indices satisfy $\ell+1\le k'+1\le n-2<n$. Hence, $\ell\in \cR(e)$. Therefore, $\srpt(e)=\min\{\ell,k'\}=\ell=\last(e)$.
\end{proof}

\begin{corollary}\label{co:ineq}
	For $e\in \bI_n(0012)$,
	$$0\le \srpt(e)\le \last(e)\le n-1.$$
\end{corollary}

\begin{proof}
	If $e= 01\cdots (n-1)$, the above inequalities are trivial since $\srpt(e)=\last(e)=n-1$. If $e\ne 01\cdots (n-1)$, the inequalities are direct consequences of Lemma \ref{le:combin-1} and the fact that $\srpt(e)\ge 0$ and $\last(e)\le n-1$.
\end{proof}

\begin{lemma}\label{le:combin-2}
	For $n\ge 2$ and $e=e_1e_2\cdots e_n\in \bI_n(0012)$, let $e$ be such that $\srpt(e)=\last(e)=k$ with $0\le k\le n-2$. Then
	
	\begin{enumerate}[leftmargin=*,align=left,labelsep=0cm,label=\textup{(\alph*)}.]
		\item for $1\le i\le k+1$,
		$$e_i=i-1;$$
		
		\item if we denote $e'=e'_1e'_2\cdots e'_{n-k}$ by the sequence obtained via $e'_i = e_{k+i}-k$ for each $1\le i\le n-k$, then $e'\in \bI_{n-k}(0012)$ such that
		$$\srpt(e')=\last(e')=0.$$
	\end{enumerate}
\end{lemma}

\begin{proof}
	Part (a) simply comes from Lemma \ref{le:combin-0}. Also, we know from Part (a) that for $k+1\le i\le n$, it holds that $e_i\ge k$. On the other hand, $e_i\le i-1$. Hence, $e'$ is still an inversion sequence. Further, it is trivial to see that $e'$ still avoids the pattern $0012$. Finally, we have $e'_1=e_{k+1}-k=k-k=0$ and $\last(e')=e'_{n-k}=e_n-k=k-k=0$. Since $n-k\ge 2>1$, we have $0\in \cR(e')$ and hence $\srpt(e')=0$.
\end{proof}

\section{Recurrences}

Let
\begin{equation*}
f_n(k,\ell):=
\left\{\begin{gathered}
\text{the number of sequences $e\in \bI_n(0012)$ with}\\[-2pt]
\text{$\srpt(e)=k$ and $\last(e)=\ell$}
\end{gathered}\right\}.
\end{equation*}
We will establish the following recurrences.

\begin{lemma}\label{le:rec}
	We have
	
	\begin{enumerate}[leftmargin=*,align=left,labelsep=0cm,label=\textup{(\alph*)}.]
		\item for $n\ge 1$,
		\begin{equation*}
		f_n(n-1,n-1)=1;
		\end{equation*}
		
		\item for $n\ge 2$,
		\begin{equation*}
		f_n(n-2,n-1)=0;
		\end{equation*}
		
		\item for $n\ge 2$ and $0\le k\le n-3$,
		\begin{equation*}
		f_n(k,n-1)=\sum_{k'=k}^{n-2}f_{n-1}(k',n-2);
		\end{equation*}
		
		\item for $n\ge 2$ and $0\le \ell\le n-2$,
		\begin{equation*}
		f_n(\ell,\ell)=\sum_{\ell'=\ell}^{n-2}\sum_{k'=\ell}^{\ell'}f_{n-1}(k',\ell');
		\end{equation*}
		
		\item for $n\ge 2$ and $0\le k<\ell\le n-2$,
		\begin{equation*}
		f_n(k,\ell)=\sum_{k'=k}^\ell f_{n-1}(k',\ell)+\sum_{\ell'=\ell}^{n-2} f_{n-1}(k,\ell').
		\end{equation*}
	\end{enumerate}
\end{lemma}

\begin{proof}
	Cases (a) and (b) are trivial. In particular, Case (a) enumerates the only inversion sequence $01\cdots (n-1)$ in which none of the letters repeat. Below we always assume that $e=e_1e_2\cdots e_n\in \bI_n(0012)$. Let $\gamma(e)$ be as in Lemma \ref{le:combin-1}.
	
	For Case (c), let $e$ be such that $\srpt(e)=k\le n-3$ and $\last(e)=n-1$. We first notice that $e_{n-1}=\last(\gamma(e))\ge \srpt(\gamma(e))$ by Corollary \ref{co:ineq}. Also, it is easy to see that $\srpt(\gamma(e))=\srpt(e)=k$ since $\last(e)=n-1>k$. Now we claim that $e_{n-1}=k$. Otherwise, namely, if $e_{n-1}>k$, we may find $i<j<n-1$ such that $e_i=e_j=k$. Hence, $e_ie_je_{n-1}e_n$ has the reduction $0012$, which contradicts the assumption that $e\in \bI_n(0012)$. We therefore have a bijection
	$$e=e_1e_2\cdots e_{n-2}(k)(n-1) \longleftrightarrow e_1e_2\cdots e_{n-2}(n-2)=e'.$$
	Notice that $e'$ is still an inversion sequence avoiding the pattern $0012$. Also, $\srpt(e')\ge k$. Otherwise, there exists some $k'<k$ that appears more than once among $e_1$, $e_2$, \ldots, $e_{n-2}$ and therefore $\srpt(e)<k$, which leads to a contradiction. Finally, to prove Case (c), it suffices to show that $e'$ could be any inversion sequence in $\bI_{n-1}(0012)$ with $\last(e')=n-2$ (which is of course true) and $\srpt(e')\ge k$. Let $e'$ be such a sequence and assume that $\srpt(e')=k'\ge k$. By Lemma \ref{le:combin-0}, we have $e_{k+1}=k$. Pulling back to $e$, we have $e_{k+1}=e_{n-1}=k$ with the indices $k+1\le n-2<n-1$. Therefore, for this $e$, we have $k\in\cR(e)$ and hence $\srpt(e)=\min\{k',k\}=k$.
	
	For Case (d), let $e$ be such that $\srpt(e)=\last(e)=\ell$ with $0\le \ell\le n-2$. We first find that $\srpt(\gamma(e))\ge \srpt(e)=\ell$. On the other hand, let $e'=e'_1e'_2\cdots e'_{n-1}\in \bI_{n-1}(0012)$ be such that $\srpt(e')\ge \ell$. By Lemma \ref{le:combin-0}, $e'_{\ell+1}=\ell$. Hence, by appending $\ell$ to the end of $e'$, we obtain a sequence with both $\srpt$ and $\last$ equal to $\ell$. We therefore arrive at a bijection between $e$ and $e'$,
	$$e=e_1e_2\cdots e_{n-1}(\ell) \longleftrightarrow e_1e_2\cdots e_{n-1}=e',$$
	and the desired relation follows.
	
	For Case (e), let $e$ be such that $\srpt(e)=k$ and $\last(e)=\ell$ with $0\le k<\ell\le n-2$. Notice that $e_{n-1}\ge k$. Otherwise, we assume that $e_{n-1}=k'<k$. Then by Lemma \ref{le:combin-0}, $e_{k'+1}=k'=e_{n-1}$. However, $k'+1<k+1<n-1$ and hence $k'\in \cR(e)$. But this violates the fact that $k=\min \cR(e)$. Now we have two cases.
	\begin{itemize}[leftmargin=*,align=left]
		\renewcommand{\labelitemi}{\scriptsize$\blacktriangleright$}
		\item $e_{n-1}<e_{n}$. We claim that $e_{n-1}=k$. Otherwise, we may find $i<j<n-1$ such that $e_i=e_j=k$. Hence, $e_ie_je_{n-1}e_n$ has the reduction $0012$, which violates the assumption that $e\in \bI_n(0012)$. Now we have a bijection between $e$ and $e'\in \bI_{n-1}(0012)$ such that $\srpt(e')\ge k$ and $\last(e')=\ell$ by
		$$e=e_1e_2\cdots e_{n-2}(k)(\ell) \longleftrightarrow e_1e_2\cdots e_{n-2}(\ell)=e'.$$
		The argument is similar to that for Case (c). This bijection leads to the first term in the right-hand side of the recurrence relation in Case (e).
		
		\item $e_{n-1}\ge e_{n}$. We have a bijection between $e$ and $e'\in \bI_{n-1}(0012)$ such that $\srpt(e')= k$ and $\last(e')\ge \ell$ by
		$$e=e_1e_2\cdots e_{n-1}(\ell) \longleftrightarrow e_1e_2\cdots e_{n-1}=e'.$$
		The argument is similar to that for Case (d). This bijection leads to the second term in the right-hand side of the recurrence relation in Case (e).
	\end{itemize}

	The proof of the lemma is therefore complete.
\end{proof}

We may therefore determine the support of $f_n(k,\ell)$.

\begin{corollary}
	For $n\ge 1$, $f_n(k,\ell)$ is supported on
	$$\{(k,\ell)\in\mathbb{N}^2:0\le k\le \ell\le n-1\}\backslash\{(n-2,n-1)\}.$$
\end{corollary}

\begin{proof}
	By Corollary \ref{co:ineq}, $f_n(k,\ell)=0$ if
	$$(k,\ell)\not\in \{(k,\ell)\in\mathbb{N}^2:0\le k\le \ell\le n-1\}.$$
	Also, $f_n(n-2,n-1)=0$ by Lemma \ref{le:rec}(b). Finally, for the remaining $(k,\ell)$, we have $f_n(k,\ell)\neq 0$ with the help of the recurrences in Lemma \ref{le:rec}.
\end{proof}

Finally, we have another recurrence.

\begin{lemma}\label{le:rec-2}
	We have, for $n\ge 2$ and $0\le k\le n-2$,
	\begin{equation*}
	f_n(k,k)=f_{n-k}(0,0).
	\end{equation*}
\end{lemma}

\begin{proof}
	This is an immediate consequence of Lemma \ref{le:combin-2}.
\end{proof}

In the sequel, we require three auxiliary functions with \textit{$q$ within a sufficiently small neighborhood of $0$}:
\begin{align*}
\cL(x;q)&:=\sum_{n\ge 1}\left(\sum_{k=0}^{n-1}f_n(k,n-1)x^k\right)q^n,\\
\cD(x;q)&:=\sum_{n\ge 1}\left(\sum_{\ell=0}^{n-2} f_n(\ell,\ell)x^\ell\right)q^n,\\
\cF(x,y;q)&:=\sum_{n\ge 1}\left(\sum_{\ell=0}^{n-1}\sum_{k=0}^{\ell} f_n(k,\ell)x^k y^\ell\right)q^n.
\end{align*}
In particular, we write, for $n\ge 1$,
\begin{align*}
L_n(x)&:=\sum_{k=0}^{n-1}f_n(k,n-1)x^k,\\
D_n(x)&:=\sum_{\ell=0}^{n-2} f_n(\ell,\ell)x^\ell,\\
F_n(x,y)&:=\sum_{\ell=0}^{n-1}\sum_{k=0}^{\ell} f_n(k,\ell)x^k y^\ell.
\end{align*}
Notice that $L_1(x)=1$, $D_1(x)=0$ and $F_1(x,y)=1$. Also, since $f_n(n-1,n-1)=1$, we have
\begin{align*}
\sum_{\ell=0}^{n-1} f_n(\ell,\ell)x^\ell = D_n(x) + x^{n-1}.
\end{align*}

\section{Proof of Theorem \ref{th:0012-last}}

Notice that Theorem \ref{th:0012-last} is equivalent to
\begin{align*}
\cL(1;q)&=\sum_{n\ge 1}\left(\sum_{k=0}^{n-1}f_n(k,n-1)\right)q^n\\
&\stackrel{?}{=} q+q^2+2q^3+4q^4+8q^5+16q^6+\cdots\\
&=\frac{q(1-q)}{1-2q}.
\end{align*}
We prove a strengthening of the above.

\begin{theorem}
	We have
	\begin{equation}\label{eq:G-gf}
	\cL(x;q)=\frac{q(1-q)^2}{(1-2q)(1-xq)}.
	\end{equation}
\end{theorem}

\begin{proof}
	For $n\ge 2$, it follows from (a), (b) and (c) of Lemma \ref{le:rec} that
	\begin{align*}
	\sum_{k=0}^{n-1}f_n(k,n-1)x^k &= x^{n-1} + \sum_{k=0}^{n-3}\sum_{k'=k}^{n-2}f_{n-1}(k',n-2)x^k\\
	&= x^{n-1} + \sum_{k'=0}^{n-3}f_{n-1}(k',n-2)\sum_{k=0}^{k'}x^k + f_{n-1}(n-2,n-2)\sum_{k=0}^{n-3}x^k\\
	&= x^{n-1} + \sum_{k'=0}^{n-3}f_{n-1}(k',n-2)\frac{1-x^{k'+1}}{1-x}+\frac{1-x^{n-2}}{1-x}.
	\end{align*}
	Therefore,
	\begin{align*}
	L_n(x)=x^{n-1}+\frac{1}{1-x}\big(L_{n-1}(1)-xL_{n-1}(x)\big)-\frac{1-x^{n-1}}{1-x}+\frac{1-x^{n-2}}{1-x}.
	\end{align*}
	Multiplying the above by $q^n$ and summing over $n\ge 2$, we have
	\begin{align*}
	\cL(x;q)-q = \frac{q}{1-x}\cL(1;q)-\frac{xq}{1-x}\cL(x;q)-\frac{q^2(1-x)}{1-xq},
	\end{align*}
	or
	\begin{align}\label{eq:func-G}
	(1-xq)(1-x+xq)\cL(x;q)=q(1-xq)\cL(1;q)+q(1-q)(1-x).
	\end{align}
	Applying the kernel method (see \cite[Exercise 4, \S{}2.2.1, p.~243]{Knu1997} or \cite{Pro2003}) yields
	\begin{align*}
	\begin{cases}
	1-x+xq=0,\\
	q(1-xq)\cL(1;q)+q(1-q)(1-x)=0.
	\end{cases}
	\end{align*}
	Solving the first equation of the system for $x$ gives
	$$x=\frac{1}{1-q}.$$
	Substituting the above into the second equation of the system, we have
	$$\cL(1;q)=\frac{q(1-q)}{1-2q}.$$
	Substituting the above back to \eqref{eq:func-G}, we arrive at \eqref{eq:G-gf}.
\end{proof}

\section{Proof of Theorem \ref{th:0012}}

We first establish two relations concerning $\cD(x;q)$.

\begin{lemma}
	We have
	\begin{align}
	\cD(x;q)&=\frac{1}{1-xq}\cD(0;q)\label{eq:D-1}\\
	&=\frac{q}{1-xq}\cF(1,1;q).\label{eq:D-2}
	\end{align}
\end{lemma}

\begin{proof}
	We know from Lemma \ref{le:rec-2} that
	\begin{align*}
	\sum_{n\ge 2}\sum_{k=0}^{n-2}f_n(k,k)x^kq^n &= \sum_{n\ge 2}\sum_{k=0}^{n-2}f_{n-k}(0,0)x^kq^n\\
	{\text{\tiny (with $n'=n-k$)}}&=\sum_{n'\ge 2}\sum_{n\ge n'}f_{n'}(0,0)x^{n-n'}q^n\\
	&=\sum_{n'\ge 2}f_{n'}(0,0) x^{-n'}\sum_{n\ge n'} (xq)^n\\
	&=\frac{1}{1-xq}\sum_{n'\ge 2}f_{n'}(0,0)q^{n'}.
	\end{align*}
	Noticing that $D_1(x)=0$, we have
	$$\cD(x;q)=\frac{1}{1-xq}\cD(0;q),$$
	which is the first part of the lemma. For the second part, we deduce from Lemma \ref{le:rec}(d) that
	\begin{align*}
	\cD(0;q)&=\sum_{n\ge 2}f_{n}(0,0)q^{n}\\
	&=\sum_{n\ge 2}\sum_{\ell'=0}^{n-2}\sum_{k'=0}^{\ell'}f_{n-1}(k',\ell')q^n\\
	&=q\cF(1,1;q).
	\end{align*}
	Therefore, \eqref{eq:D-2} follows.
\end{proof}

Next, we show a relation between $\cF(x,1;q)$ and $\cF(1,1;q)$.

\begin{lemma}
	We have
	\begin{equation}\label{eq:F(x,1;q)}
	\cF(x,1;q)=\frac{1-q}{1-xq}\cF(1,1;q).
	\end{equation}
\end{lemma}

\begin{proof}
	For $n\ge 2$, it follows from Lemma \ref{le:rec}(d) that
	\begin{align*}
	D_n(x)&=\sum_{\ell=0}^{n-2} f_n(\ell,\ell)x^\ell\notag\\
	&= \sum_{\ell=0}^{n-2} \sum_{\ell'=\ell}^{n-2}\sum_{k'=\ell}^{\ell'}f_{n-1}(k',\ell')x^\ell\notag\\
	&=\sum_{\ell'=0}^{n-2}\sum_{k'=0}^{\ell'}f_{n-1}(k',\ell') \sum_{k=0}^{k'}
	x^\ell\notag\\
	&=\sum_{\ell'=0}^{n-2}\sum_{k'=0}^{\ell'}f_{n-1}(k',\ell') \frac{1-x^{k'+1}}{1-x}\notag\\
	&=\frac{1}{1-x}\big(F_{n-1}(1,1)-xF_{n-1}(x,1)\big).
	\end{align*}
	Therefore,
	$$\cD(x;q)=\frac{q}{1-x}\big(\cF(1,1;q)-x\cF(x,1;q)\big).$$
	Substituting \eqref{eq:D-2} into the above yields
	$$\frac{q}{1-xq}\cF(1,1;q)=\frac{q}{1-x}\big(\cF(1,1;q)-x\cF(x,1;q)\big),$$
	from which \eqref{eq:F(x,1;q)} follows.
\end{proof}

We then construct a functional equation for $\cF(x,y;q)$.

\begin{lemma}
	We have
	\begin{align}\label{eq:F(x,y;q)}
	&\left(1+\frac{xq}{1-x}+\frac{yq}{1-y}\right)\cF(x,y;q)\notag\\
	&\quad= \frac{q}{1-x}\cF(1,y;q)+\frac{q(1-q)}{(1-y)(1-xyq)}\cF(1,1;q)+\frac{q(1-q-2yq+2yq^2 +y^2q^2)}{(1-2yq)(1-xyq)}.
	\end{align}
\end{lemma}

\begin{proof}
	We first observe that
	\begin{align}\label{eq:lemma-4.3-1}
	\sum_{\ell=0}^{n-2} f_n(\ell,\ell)x^\ell y^\ell+\sum_{\ell=1}^{n-2}\sum_{k=0}^{\ell-1}f_n(k,\ell)x^k y^\ell &= F_n(x,y)-\sum_{k=0}^{n-1}f_n(k,n-1)x^k y^{n-1}\notag\\
	&=F_n(x,y)- y^{n-1}L_n(x).
	\end{align}
	Notice also that
	\begin{equation}\label{eq:lemma-4.3-2}
	\sum_{\ell=0}^{n-2} f_n(\ell,\ell)x^\ell y^\ell = D_n(xy).
	\end{equation}
	
	Now, by Lemma \ref{le:rec}(e), we may separate
	\begin{align*}
	\sum_{\ell=1}^{n-2}\sum_{k=0}^{\ell-1}f_n(k,\ell)x^k y^\ell &= \sum_{\ell=1}^{n-2}\sum_{k=0}^{\ell-1}\sum_{k'=k}^\ell f_{n-1}(k',\ell)x^k y^\ell\\
	&\quad+\sum_{\ell=1}^{n-2}\sum_{k=0}^{\ell-1}\sum_{\ell'=\ell}^{n-2} f_{n-1}(k,\ell')x^k y^\ell.
	\end{align*}
	We further notice that the first term on the right-hand side can be separated as
	\begin{align*}
	\sum_{\ell=1}^{n-2}\sum_{k=0}^{\ell-1}\sum_{k'=k}^\ell f_{n-1}(k',\ell)x^k y^\ell=\sum_{\ell=1}^{n-2}\sum_{k=0}^{\ell-1}\sum_{k'=k}^{\ell-1} f_{n-1}(k',\ell)x^k y^\ell+\sum_{\ell=1}^{n-2}\sum_{k=0}^{\ell-1}f_{n-1}(\ell,\ell)x^ky^\ell.
	\end{align*}
	We have
	\begin{align*}
	&\sum_{\ell=1}^{n-2}\sum_{k=0}^{\ell-1}\sum_{k'=k}^{\ell-1} f_{n-1}(k',\ell)x^k y^\ell\\
	&\quad=\sum_{\ell=1}^{n-2}\sum_{k'=0}^{\ell-1} f_{n-1}(k',\ell) y^\ell \sum_{k=0}^{k'}x^k\\
	&\quad=\sum_{\ell=1}^{n-2}\sum_{k'=0}^{\ell-1} f_{n-1}(k',\ell) y^\ell \frac{1-x^{k'+1}}{1-x}\\
	&\quad=\sum_{\ell=0}^{n-2}\sum_{k'=0}^{\ell} f_{n-1}(k',\ell) y^\ell \frac{1-x^{k'+1}}{1-x}-\sum_{\ell=0}^{n-2}f_{n-1}(\ell,\ell)y^\ell\frac{1-x^{\ell+1}}{1-x}\\
	&\quad=\frac{1}{1-x}\big(F_{n-1}(1,y)-xF_{n-1}(x,y)\big)\\
	&\quad\quad-\frac{1}{1-x}\big(D_{n-1}(y)+y^{n-2}-xD_{n-1}(xy)-x^{n-1}y^{n-2}\big).
	\end{align*}
	Also,
	\begin{align*}
	\sum_{\ell=1}^{n-2}\sum_{k=0}^{\ell-1}f_{n-1}(\ell,\ell)x^ky^\ell&=\sum_{\ell=1}^{n-2}f_{n-1}(\ell,\ell)y^\ell\frac{1-x^\ell}{1-x}\\
	&=\sum_{\ell=0}^{n-2}f_{n-1}(\ell,\ell)y^\ell\frac{1-x^\ell}{1-x}\\
	&=\frac{1}{1-x}\big(D_{n-1}(y)+y^{n-2}-D_{n-1}(xy)-x^{n-2}y^{n-2}\big).
	\end{align*}
	On the other hand,
	\begin{align*}
	\sum_{\ell=1}^{n-2}\sum_{k=0}^{\ell-1}\sum_{\ell'=\ell}^{n-2} f_{n-1}(k,\ell')x^k y^\ell&=\sum_{\ell'=1}^{n-2}\sum_{k=0}^{\ell'-1}f_{n-1}(k,\ell')x^k \sum_{\ell=k+1}^{\ell'}y^\ell\\
	&=\sum_{\ell'=1}^{n-2}\sum_{k=0}^{\ell'-1}f_{n-1}(k,\ell')x^k\frac{y^{k+1}-y^{\ell'+1}}{1-y}\\
	&=\sum_{\ell'=0}^{n-2}\sum_{k=0}^{\ell'}f_{n-1}(k,\ell')x^k\frac{y^{k+1}-y^{\ell'+1}}{1-y}\\
	&=\frac{y}{1-y}\big(F_{n-1}(xy,1)-F_{n-1}(x,y)\big).
	\end{align*}
	Therefore,
	\begin{align}\label{eq:lemma-4.3-3}
	&\sum_{\ell=1}^{n-2}\sum_{k=0}^{\ell-1}f_n(k,\ell)x^k y^\ell \notag\\
	&\quad=\frac{1}{1-x}\big(F_{n-1}(1,y)-xF_{n-1}(x,y)\big)+\frac{y}{1-y}\big(F_{n-1}(xy,1)-F_{n-1}(x,y)\big)\notag\\
	&\quad\quad-D_{n-1}(xy)-x^{n-2}y^{n-2}.
	\end{align}
	
	It follows from \eqref{eq:lemma-4.3-1}, \eqref{eq:lemma-4.3-2} and \eqref{eq:lemma-4.3-3} that
	\begin{align*}
	&F_n(x,y)- y^{n-1}L_n(x)\\
	&\quad=D_n(xy)+\frac{1}{1-x}\big(F_{n-1}(1,y)-xF_{n-1}(x,y)\big)\notag\\
	&\quad\quad +\frac{y}{1-y}\big(F_{n-1}(xy,1)-F_{n-1}(x,y)\big)-D_{n-1}(xy)-x^{n-2}y^{n-2}.
	\end{align*}
	Therefore,
	\begin{align*}
	&\cF(x,y;q)-y^{-1}\cL(x;yq)\\
	&\quad=\cD(xy;q)+\frac{q}{1-x}\big(\cF(1,y;q)-x\cF(x,y;q)\big)\\
	&\quad\quad +\frac{yq}{1-y}\big(\cF(xy,1;q)-\cF(x,y;q)\big)-q\cD(xy;q)-\frac{q^2}{1-xyq}.
	\end{align*}
	Applying \eqref{eq:G-gf}, \eqref{eq:D-2} and \eqref{eq:F(x,1;q)} gives the desired result.
\end{proof}

With the assistance of the kernel method, we may deduce a functional equation satisfied by $\cF(1,y;q)$.

\begin{lemma}
	We have
	\begin{align}\label{eq:F(1,y;q)}
	\cF(1,y;q)&=\frac{q}{1-y+y^2 q}\cF(1,1;q)+\frac{q(1-y)(1-q-2yq+2yq^2 +y^2q^2)}{(1-q)(1-2yq)(1-y+y^2 q)}.
	\end{align}
\end{lemma}

\begin{proof}
	We multiply both sides of \eqref{eq:F(x,y;q)} by $(1-x)(1-y)$. Then
	\begin{align*}
	&\big((1-y+yq)-x(1-y-q+2yq)\big)\cF(x,y;q)\\
	&\quad= q(1-y)\cF(1,y;q)+\frac{q(1-q)(1-x)}{1-xyq}\cF(1,1;q)\\
	&\quad\quad+\frac{q(1-x)(1-y)(1-q-2yq+2yq^2 +y^2q^2)}{(1-2yq)(1-xyq)}.
	\end{align*}
	We treat the kernel polynomial as a function in $x$ and solve
	$$(1-y+yq)-x(1-y-q+2yq)=0$$
	so that
	$$x=\frac{1-y+yq}{1-y-q+2yq}.$$
	Substituting the above into
	\begin{align*}
	0&= q(1-y)\cF(1,y;q)+\frac{q(1-q)(1-x)}{1-xyq}\cF(1,1;q)\\
	&\quad+\frac{q(1-x)(1-y)(1-q-2yq+2yq^2 +y^2q^2)}{(1-2yq)(1-xyq)},
	\end{align*}
	we arrive at \eqref{eq:F(1,y;q)} after simplification.
\end{proof}

Finally, we are ready to complete the proof of Theorem \ref{th:0012}.

\begin{proof}[Proof of Theorem \ref{th:0012}]
	It is known that (cf.~\cite[A279561]{OEIS})
	\begin{align}
	1+\sum_{n\ge 1}\Bigg(1+\sum_{i=1}^{n-1}\binom{2i}{i-1}\Bigg)q^n=\frac{1-4q+(1-2q)\sqrt{1-4q}}{2(1-q)(1-4q)}.
	\end{align}
	We then rewrite \eqref{eq:F(1,y;q)} as
	\begin{align*}
	(1-y+y^2 q)\cF(1,y;q)=q\cF(1,1;q)+\frac{q(1-y)(1-q-2yq+2yq^2 +y^2q^2)}{(1-q)(1-2yq)}.
	\end{align*}
	We treat the kernel polynomial as a function in $y$ and solve
	$$1-y+y^2 q = 0.$$
	Then
	$$y_{1,2}=\frac{1\mp \sqrt{1-4q}}{2q}.$$
	We choose the solution
	$$y_{1}=\frac{1- \sqrt{1-4q}}{2q}$$
	since $y_1\to 0$ as $q\to 0$. Substituting $y=y_1$ into
	$$0=q\cF(1,1;q)+\frac{q(1-y)(1-q-2yq+2yq^2 +y^2q^2)}{(1-q)(1-2yq)},$$
	we find that
	\begin{align}
	\cF(1,1;q)&=\frac{-(1-2q)(1-4q)+(1-2q)\sqrt{1-4q}}{2(1-q)(1-4q)}\notag\\
	&=\frac{1-4q+(1-2q)\sqrt{1-4q}}{2(1-q)(1-4q)}-1.
	\end{align}
	This implies that for $n\ge 1$,
	\begin{align*}
	1+\sum_{i=1}^{n-1}\binom{2i}{i-1}=\sum_{\ell=0}^{n-1}\sum_{k=0}^{\ell} f_n(k,\ell)=|\bI_n(0012)|.
	\end{align*}
	Therefore, Conjecture \ref{conj:Lin-Ma} is true.
\end{proof}

\bibliographystyle{amsplain}

\end{document}